\documentclass[12pt]{amsart}

\usepackage{amsmath,amssymb,amscd,amsthm}
\usepackage{hyperref}

\makeatletter
\@namedef{subjclassname@2010}{%
  \textup{2010} Mathematics Subject Classification}
\makeatother

\newtheorem{theorem}{Theorem}[section]

\newtheorem {corollary}[theorem]{Corollary}

\theoremstyle{definition}

\newtheorem{example}{Example}[section]

\DeclareMathOperator\Sym{Sym}

\newcommand{\PP}{\mathbb P}

\newcommand{\NN}{\mathbb N}

\begin{document}

\baselineskip=16pt

\title[The degree of Fano schemes]{On the degree of Fano schemes\\ of linear subspaces on hypersurfaces}

\author{Dang Tuan Hiep}

\address{Faculty of Mathematics and Computer Science, University of Dalat, 01 Phu Dong Thien Vuong, Da Lat, Vietnam}

\email{hiepdt@dlu.edu.vn}

\address{National Center for Theoretical Sciences, Mathematics Division, 2F of Astronomy-Mathematics Building, National Taiwan University, Taipei 10617, Taiwan}

\begin{abstract}
In this paper we propose and prove an explicit formula for computing the degree of Fano schemes of linear subspaces on general hypersurfaces. The method used here is based on the localization theorem and Bott's residue formula in equivariant intersection theory.
\end{abstract}

\subjclass[2010]{14C15; 14N15; 55N91}

\keywords{Fano schemes, Schubert calculus, equivariant intersection theory, Bott's residue formula}

\date{\today}

\maketitle

\section{Introduction}

Let $X \subset \mathbb P_{\mathbb C}^n$ be a general hypersurface of degree $d$. The {\it Fano scheme} $F_k(X)$ is defined to be the set of $k$-dimensional linear subspaces of $\mathbb P_{\mathbb C}^n$ which are contained in $X$. This is a subscheme of the Grassmannian $\mathbb G(k,n)$ of $k$-dimensional linear subspaces in $\mathbb P_{\mathbb C}^n$. For convenience, we set
$$\delta = (k+1)(n-k)-\binom{d+k}{d}.$$
Suppose that $d \neq 2$ (or $n \geq 2k +1$) and $\delta \geq 0$. Langer \cite{L} showed that $F_k(X)$ is smooth of expected dimension $\delta$. In this case, the degree of $F_k(X)$ is given by the following formula
\begin{equation}\label{degreeFano}
\deg(F_k(X)) = \int_{\mathbb G(k,n)}c_{\binom{d+k}{d}}(\Sym^dS^\vee)\cdot c_1(Q)^{\delta},
\end{equation}
where $S$ and $Q$ are respectively the tautological sub-bundle and quotient bundle on $\mathbb G(k,n)$, $\Sym^dS^\vee$ is the $d$-th symmetric power of the dual of $S$, and $c_i(E)$ is the $i$-th Chern class of the vector bundle $E$. Note that $\int_Y\alpha$ denotes the degree of the cycle class $\alpha$ on $Y$ defined in \cite[Definition 1.4]{F}. Formula \eqref{degreeFano} can be found, for example, in \cite[Example 14.7.13]{F} or \cite[Section 3.5]{M}. Using Schubert calculus, Debarre and Manivel \cite[Theorem 4.3]{DM} showed that the degree of $F_k(X)$ is equal to a certain coefficient of an explicit polynomial, given as the product of linear forms. In this paper we propose and prove an explicit formula for computing the degree of $F_k(X)$ via equivariant intersection theory. 

We denote by $\mathcal I$ the set of $(k+1)$-subsets of $\{1,\ldots,n+1\}$, where a $(k+1)$-subset of $\{1,\ldots,n+1\}$ is a subset consisting of $k+1$ elements of $\{1,\ldots,n+1\}$. Consider the polynomial ring $\mathbb C[h_1,\ldots,h_{n+1}]$ in $n+1$ variables $h_1,\ldots,h_{n+1}$. For each $I\in \mathcal I$, we set
$$S_I = \prod_{v_i\in \NN, \sum_{i\in I} v_i = d}\left(\sum_{i\in I} v_i h_i\right), Q_I = \sum_{j\not\in I}h_j,$$
and
$$T_I = \prod_{i\in I} \prod_{j\not\in I}(h_i - h_j).$$
Then $S_I, Q_I,$ and $T_I$ are elements in $\mathbb C[h_1,\ldots,h_{n+1}]$.

Here is the main result of this paper:

\begin{theorem}\label{degreeFanoBott}
Let $k,d,n \in \NN$ satisfy $d\neq 2$ (or $n \geq 2k+1$) and $\delta \geq 0$, and let $X \subset \PP_{\mathbb C}^n$ be a general hypersurface of degree $d$. Then the degree of the Fano scheme $F_k(X)$ of $k$-dimensional linear subspaces on $X$ is given by the following formula:
$$\deg(F_k(X)) = (-1)^\delta\sum_{I \in \mathcal I}\frac{S_I Q_I^\delta}{T_I}.$$
\end{theorem}

The right-hand-side of the formula in Theorem \ref{degreeFanoBott} is a rational polynomial function, and the above theorem claims in other words that it is in fact a constant function, moreover it is an integer. Namely, for any numbers $h_i$ such that $h_i \neq h_j$ for $i \neq j$, the right-hand-side of the formula becomes the same integer. Indeed, for concrete numbers $(k,d,n)$ such that the above conditions are satisfied, it can be checked by computer algebra systems. See \cite[Chapter 5]{D} for the computation in \textsc{Sage} \cite{S} and \textsc{Singular} \cite{DGPS}.

If $k,d,n \in \mathbb{N}$ satisfy $d \neq 2$ (or $n\geq 2k+1$) and $\delta = 0$, then the Fano scheme $F_k(X)$ is zero-dimensional. In this case, the degree of $F_k(X)$ is equal to the number of $k$-dimensional linear subspaces on $X$. 

In particular, if $k=1$ and $\delta = 0$, then $d=2n-3$. In this case, we have 
$$\mathcal I = \{\{i,j\} \mid 1 \leq i < j \leq n+1\}.$$
For each $I = \{i,j\} \in \mathcal I$, we have
$$S_I = \prod_{a = 0}^{2n-3}(ah_i+(2n-3-a)h_j)$$
and
$$T_I = \prod_{k\neq i,j}(h_i-h_k)(h_j-h_k).$$
Thus we have the following corollary.
\begin{corollary}\label{linearBott}
Let $k,d,n \in \mathbb{N}$ satisfy $d \neq 2$ (or $n\geq 2k+1$) and $\delta = 0$, and let $X\subset \PP_{\mathbb C}^n$ be a general hypersurface of degree $d$. Then the number of $k$-dimensional linear subspaces on $X$ is equal to 
$$\sum_{I \in \mathcal I}\frac{S_I}{T_I}.$$
In particular, the number of lines on a general hypersurface of degree $2n-3$ in $\mathbb P^n_{\mathbb C}$ is equal to
$$\sum_{1 \leq i < j \leq n+1}\frac{\displaystyle\prod_{a = 0}^{2n-3}(ah_i+(2n-3-a)h_j)}{\displaystyle \prod_{k\neq i,j}(h_i-h_k)(h_j-h_k)}.$$
\end{corollary}

Debarre and Manivel \cite{DM} showed that the degree of $F_k(X)$ is equal to the coefficient of the monomial $x_0^nx_1^{n-1}\cdots x_k^{n-k}$ in the product of the polynomial 
$$\prod_{v_i\in \NN, \sum_{i=0}^k v_i = d}\left(\sum_{i=0}^k v_i x_i\right)(x_0+\cdots+x_k)^\delta$$
by the Vandermonde determinant 
$$\prod_{0\leq i < j \leq k}(x_i - x_j).$$
See \cite[Theorem 3.5.18]{M} for the special case of this formula for computing the number of linear subspaces on hypersurfaces. 

In \cite{AK}, Altman and Kleiman gave a formula for computing the degree of $F_1(X)$. Harris \cite{H} gave a formula for computing the number of lines on a general hypersurface of degree $2n-3$ in $\mathbb P^n$. In \cite{W}, van der Waerden showed that the number of lines on a general hypersurface of degree $d = 2n-3$ in $\mathbb P^n_{\mathbb C}$ could also be computed as the coefficient of the monomial $x^ny^{n-1}$ in the polynomial 
$$(x-y)\prod_{i=0}^d((d-i)x+iy).$$
These results were obtained via Schubert calculus. However, in this paper, our approach is completely different. We use the localization theorem and Bott's residue formula in equivariant intersection theory.

The rest of the paper is organized as follows: Section 2 presents a quick review of equivariant intersection theory. Section 3 presents the proof of the main theorem. 

\section{Equivariant intersection theory}

Edidin and Graham \cite{EG1,EG2} gave an algebraic construction to equivariant intersection theory. In this section, we review the basic notions and results of this theory. Let $G$ be a linear algebraic group and let $X$ be a scheme of finite type over $\mathbb C$ endowed with a $G$-action. For any non-negative integer $i$, we can find a representation $V$ of $G$ together with a dense open subset $U \subset V$ on which $G$ acts freely and whose complement has codimension larger than $\dim X - i$ such that the principal bundle quotient $U \rightarrow U/G$ exists in the category of schemes (see \cite[Lemma 9]{EG1}). The diagonal action on $X \times U$ is then also free, which implies that under mild assumption, a principal bundle quotient $X \times U \rightarrow (X \times U)/G$ exists in the category of schemes (see \cite[Proposition 23]{EG1}). In what follows, we will tacitly assume that the scheme $(X \times U)/G$ exists and denote it by $X_G$.
\subsection{Equivariant Chow groups}
We define the $i$-th $G$-equivariant Chow group of $X$ to be
$$A_i^G(X) := A_{i+\dim U - \dim G}(X_G),$$
where $A_*$ stands for the ordinary Chow group defined in \cite{F}. By \cite[Proposition 1]{EG1}, this is well-defined. The $G$-equivariant Chow group of $X$ is defined to be
$$A_*^G(X) = \bigoplus_i A_i^G(X).$$
If $X_G$ is smooth, then $A_*^G(X)$ inherits an intersection product from the ordinary Chow groups. This endows $A_*^G(X)$ with the structure of a graded ring, called the {\it $G$-equivariant Chow ring} of $X$. For example, if $G = T = (\mathbb C^*)^n$ is a split torus of dimension $n$, then the $T$-equivariant Chow ring of a point is isomorphic to a polynomial ring in $n$ variables (see \cite[Section 3.2]{EG1}). Throughout this paper, we denote this ring by $R_T$.
\subsection{Equivariant vector bundles and Chern classes}
 \label{equiBundle}
 
A $G$-equivariant vector bundle is a vector bundle $E$ on $X$ such that the action of $G$ on $X$ lifts to an action of $G$ on $E$ which is linear on fibers. By \cite[Lemma 1]{EG1}, $E_G$ is a vector bundle over $X_G$. The {\it $G$-equivariant Chern classes} $c_i^G(E)$ are defined to be the Chern classes $c_i(E_G)$. If $E$ has rank $r$, then the top Chern class $c_r^G(E)$ is called the {\it $G$-equivariant Euler class} of $E$ and denoted $e^G(E)$.

Note that a $G$-equivariant vector bundle over a point is a representation of $G$ (see also \cite[Section 3.2]{EG1}). Our primary interest is when $G = T = (\mathbb C^*)^n$ is a split torus of dimension $n$ and $X = pt$ is a point. In this case, let $M(T)$ be the character group of the torus $T$. Suppose that $R_T = \mathbb C[h_1,\ldots,h_n]$. There is a group homomorphism $\psi : M(T) \rightarrow R_T$ given by $\rho_i \mapsto h_i,$ where $\rho_i$ is the character of $T$ defined by $\rho_i(t_1,\ldots,t_n) = t_i$. This induces a ring isomorphism $\Sym(M(T)) \simeq R_T$. We call $\psi(\rho)$ the {\it weight} of $\rho$. In particular, $h_i$ is the weight of $\rho_i$.

\begin{example}\cite[Example 9.1.1.1]{CK}
The diagonal action of $G = T = (\mathbb C^*)^n$ on $\mathbb C^n$ gives a $T$-equivariant vector bundle $E$ over a point. The corresponding representation of $T$ has characters $\rho_i$ for $i = 1, \ldots,n$, and their weights are $h_i$. In this case, we have $E_T \cong \mathcal O(h_1) \oplus \cdots \oplus \mathcal O(h_n).$ This implies that
$c_i^T(E) = c_i(E_T) = s_i(h_1, \ldots, h_n) \in R_T,$ where $s_i$ is the $i$-th elementary symmetric function. In particular, the $T$-equivariant Euler class of $E$ is
$e^T(E) = h_1 \cdots h_n \in R_T.$
\end{example}

\begin{example}
Consider the diagonal action of $T = (\mathbb C^*)^4$ on $\mathbb C^4$ given in coordinates by
$$(t_1,t_2,t_3,t_4) \cdot (x_1, x_2 , x_3, x_4) = (t_1x_1, t_2x_2, t_3x_3, t_4x_4).$$
This induces an action of $T$ on the Grassmannian $G(2,4)$ with the isolated fixed points $L_I$ corresponding to coordinate $2$-planes in $\mathbb C^4$. Each $L_I$ is indexed by the $2$-subset $I$ of the set $\{1,2,3,4\}$ so that $L_I$ is defined by the equations $x_j=0$ for $j\not \in I$. Let $S$ be the tautological sub-bundle on $G(2,4)$. At each $L_I$, the restriction of the action of $T$ on the fiber $S|_{L_I}$ gives a representation of $T$ with characters $\rho_i$ for $i \in I$. This representation gives a $T$-equivariant vector bundle of rank $2$ over a point. We also denote it by $S|_{L_I}$. If $I = \{i_1,i_2\}$, then we have $c_1^T(S|_{L_I}) = h_{i_1} + h_{i_2}$ and $c_2^T(S|_{L_I}) = h_{i_1} \cdot h_{i_2},$ where $h_{i_1}$ and $h_{i_2}$ are the weights of $\rho_{i_1}$ and $\rho_{i_2}$ respectively.
\end{example}
\subsection{Localization and Bott's residue formula}
Let $X$ be a scheme endowed with an action of $T = (\mathbb C^*)^n$. We denote the fixed point locus by $X^T$. The localization theorem states that up to $R_T$-torsion, the $T$-equivariant Chow group of the fixed points locus $X^T$ is isomorphic to that of $X$. Moreover, the localization isomorphism is given by the equivariant push-forward induced by the inclusion of $X^T$ to $X$ (see \cite[Theorem 1]{EG2}). For smooth varieties, the inverse to the equivariant push-forward can be written explicitly (see \cite[Theorem 2]{EG2}). Using these results, Edidin and Graham gave an algebraic proof of Bott's residue formula for Chern numbers of vector bundles on smooth complete varieties (see \cite[Theorem 3]{EG2}). Bott's residue formula shows that we can compute the degree of a zero-dimensional cycle class on a smooth complete variety $X$ in terms of local contributions supported on the components of the fixed point locus of a torus action on $X$. 

Let us describe Bott's residue formula in case where the fixed point locus $X^T$ is finite. Let $E$ be a $T$-equivariant vector bundle on $X$. For each point $Z_j \in X^T$, assume that the action of $T$ on $E|_{Z_j}$ has characters $\chi_j^1,\ldots,\chi_j^s$, where each $\chi_j^i$ is a linear combination of the basic characters $\rho_1,\ldots,\rho_n$ defined above. Write $\chi_j^i = a_j^i(\rho_1,\ldots,\rho_n)$ for this linear combination. Since $\rho_i$ has weight $h_i$, the weight of $\chi_j^i$ is $a_j^i(\bar h)$, where $a_j^i(\bar h) \in R_T$ denotes the linear combination of $h_1,\ldots,h_n$ obtained by replacing $h_i$ by $\rho_i$ in $a_j^i(\rho_1,\ldots,\rho_n)$. Thus the $T$-action on $E|_{Z_j}$ has weights $a_j^i(\bar h)$. This implies that
$$c_k^T(E|_{Z_j}) = s_k(a_j^1(\bar h), \ldots, a_j^s(\bar h)) \in R_T,$$
where $s_k$ is the $k$-th elementary symmetric function.
In addition, since $Z_j$ is just a point, the normal bundle $N_j$ of $Z_j$ in $X$ is the tangent space to $X$ at $Z_j$. If the action of $T$ on the tangent space to $X$ at $Z_j$ has weights $b_j^1(\bar h), \ldots, b_j^m(\bar h)$, where $m = \dim X$, then
$$e^T(N_j) = \prod_{i=1}^m b_j^i(\bar h) \in R_T.$$
We arrive at the following formula
\begin{equation}\label{bott}
\int_X P(c_k(E)) = \sum_j \frac{P(s_k(a_j^1(\bar h), \ldots, a_j^s(\bar h)))}{\displaystyle \prod_{i=1}^m b_j^i(\bar h)},
\end{equation}
where $P(c_k(E))$ is a polynomial in the Chern classes of $E$.

Let us look closer at a simple example based on Bott's residue formula for the projective space $\mathbb P_{\mathbb C}^2$.

\begin{example}
Consider the natural action of $T = (\mathbb C^*)^3$ on $\mathbb P_{\mathbb C}^2$ given in coordinates by 
$$(t_1,t_2,t_3) \cdot (x_1: x_2 : x_3) = (t_1x_1: t_2x_2 : t_3x_3).$$
The fixed point locus is finite and consists of coordinate lines, say $p_1,p_2,$ and $p_3$. Let $c$ be the first Chern class of $E = \mathcal O_{\mathbb P_{\mathbb C}^2}(1)$, that is the class of a hyperplane. It is well-known that 
$$\int_{\mathbb P_{\mathbb C}^2}c^2 = 1.$$
Applying Bott's residue formula (\ref{bott}) we get the identity
$$\frac{h_1^2}{(h_3-h_1)(h_2-h_1)}+\frac{h_2^2}{(h_3-h_2)(h_1-h_2)}+\frac{h_3^2}{(h_1-h_3)(h_2-h_3)} = 1,$$
which can be checked by hand.
\end{example}

\section{The proof of Theorem \ref{degreeFanoBott}}

Consider the diagonal action of $T = (\mathbb C^*)^{n+1}$ on $\mathbb P_{\mathbb C}^n$ given in coordinates by
$$(t_1,\ldots,t_{n+1}) \cdot (x_1: \cdots : x_{n+1}) = (t_1x_1: \cdots : t_{n+1}x_{n+1}).$$
This induces an action of $T$ on the Grassmannian $\mathbb G(k,n)$ with $\binom{n+1}{k+1}$ isolated fixed points $L_I$ corresponding to $\binom{n+1}{k+1}$ coordinate $k$-planes in $\PP_{\mathbb C}^n$. Each fixed point $L_I$ is indexed by a $(k+1)$-subset $I$ of the set $\{1,\ldots,n+1\}$.
Let $S$ and $Q$ be the tautological sub-bundle and quotient bundle on $\mathbb G(k,n)$ respectively. The key idea is that, at each $L_I$, the torus action on the fibers $S|_{L_I}$ and $Q|_{L_I}$ have characters $\rho_i$ for $i\in I$ and $\rho_j$ for $j\not\in I$ respectively. Since the tangent bundle on the Grassmannian is isomorphic to $S^\vee \otimes Q$, the characters of the torus action on the tangent space at $L_I$ are
$$\{\rho_j - \rho_i \mid i \in I, j \not \in I \}.$$
The normal bundle $N_{L_I}$ of $L_I$ in $\mathbb G(k,n)$ is just the tangent space of $\mathbb G(k,n)$ at $L_I$. Hence
\begin{align*}
e^T(N_{L_I}) & = \prod_{i\in I} \prod_{j\not\in I} (h_j-h_i)\\
& = (-1)^{(k+1)(n-k)}\prod_{i\in I} \prod_{j\not\in I} (h_i-h_j),
\end{align*}
where $h_i$ is the weight of $\rho_i$ defined above, and $R_T = \mathbb C[h_1,\ldots,h_{n+1}]$ is the $T$-equivariant Chow ring of a point.

At each $L_I$, we also need to compute $e^T(\Sym^dS^\vee|_{L_I})$ and $c_1^T(Q|_{L_I})$. Since the characters of the torus action on $S^\vee|_{L_I}$ are $-\rho_i$ for $i\in I$, the torus action on $\Sym^dS^\vee|_{L_I}$ has characters 
$$\left\{\sum_{i\in I} v_i(-\rho_i) \mid v_i\in \NN, \sum_{i\in I} v_i = d \right\}.$$
Hence
\begin{align*}
e^T(\Sym^dS^\vee|_{L_I}) & = \prod_{v_i\in \NN, \sum_{i\in I} v_i = d}\left(\sum_{i\in I} v_i(-h_i)\right) \\
& = (-1)^{\binom{d+k}{d}}\prod_{v_i\in \NN, \sum_{i\in I} v_i = d}\left(\sum_{i\in I} v_ih_i\right).
\end{align*}
Since the characters of the torus action on $Q|_{L_I}$ are $\rho_j$ for $j \not\in I$, we have
$$c^T_1(Q|_{L_I}) = \sum_{j \not\in I}h_j.$$
By \eqref{degreeFano} and \eqref{bott}, we obtain 
\begin{align*}
\deg(F_k(X)) & = \sum_{I\in \mathcal I} \frac{e^T(\Sym^dS^\vee|_{L_I})(c^T_1(Q|_{L_I}))^\delta}{e^T(N_{L_I})} \\
& = (-1)^\delta\sum_{I \in \mathcal I}\frac{\left(\displaystyle\prod_{v_i\in \NN, \sum_{i\in I} v_i = d}\left(\sum_{i\in I} v_ih_i\right)\right)\left(\displaystyle\sum_{j \not\in I}h_j \right)^{\delta}}{\displaystyle \prod_{i\in I}\prod_{j\not\in I}(h_i-h_j)}.
\end{align*}
We get the desired formula.

Another proof of Theorem \ref{degreeFanoBott} is as follows. Consider the action of $T = \mathbb C^*$ on $\mathbb P_{\mathbb C}^n$ given in coordinates by
$$t \cdot (x_1: \cdots : x_{n+1}) = (t^{h_1}x_1: \cdots : t^{h_{n+1}}x_{n+1}).$$
The induced action of $T$ on $\mathbb G(k,n)$ also has $\binom{n+1}{k+1}$ isolated fixed points $L_I$ as above. At each $L_I$, the torus action on the fibers $S|_{L_I}$ and $Q|_{L_I}$ have characters $h_i \rho$ for $i\in I$ and $h_j \rho$ for $j\not\in I$ respectively, where $\rho$ is the character of $T$ defined by $\rho(t) = t$ for all $t \in T$. We denote the weight of $\rho$ by $h$. In this case, the $T$-equivariant Chow ring of a point is $R_T = \mathbb C[h]$. With this set-up, we have
\begin{align*}
e^T(N_{L_I}) & = \prod_{i\in I} \prod_{j\not\in I} (h_j-h_i) h^{(k+1)(n-k)}\\
& = (-1)^{(k+1)(n-k)}T_I h^{(k+1)(n-k)}.
\end{align*}
Similarly, we also have
\begin{align*}
e^T(\Sym^dS^\vee|_{L_I}) & = \prod_{v_i\in \NN, \sum_{i\in I} v_i = d}\left(\sum_{i\in I} v_i(-h_i h)\right) \\ 
& = (-1)^{\binom{d+k}{d}} S_I h^{\binom{d+k}{d}},
\end{align*}
and
$$c^T_1(Q|_{L_I}) = \sum_{j \not\in I}h_j h = Q_I h.$$
By \eqref{degreeFano} and \eqref{bott}, we obtain
$$\deg(F_k(X)) = (-1)^\delta\sum_{I \in \mathcal I}\frac{S_I Q_I^\delta h^{(k+1)(n-k)}}{T_I h^{(k+1)(n-k)}}.$$
Cancelling $h^{(k+1)(n-k)}$, we get the desired formula.

\section*{Acknowledgements}

This work is a part of my Ph.D. thesis at the University of Kaiserslautern. I would like to take this opportunity to express my profound gratitude to my advisor Professor Wolfram Decker. I would also like to thank Dr. Janko B\"{o}hm for his valuable suggestions.

\end{document}